\documentclass[letterpaper, 10 pt, conference]{ieeeconf}  

\usepackage{svg}
\usepackage{comment}
\usepackage{lipsum}
\usepackage{graphicx}
\usepackage{amsmath,amsfonts,amssymb,amscd}
\usepackage[colorlinks]{hyperref}
\usepackage{mathtools}
\usepackage{notes2bib}

\newtheorem{thm}{Theorem}
\newtheorem{lem}{Lemma}
\newtheorem{proposition}{Proposition}

\newtheorem{definition}{Definition}

\usepackage{algorithm}
\usepackage{algpseudocode}

\newcommand{\R}{\mathbb{R}}
\newcommand{\Rn}{\R^n}
\newcommand{\Rm}{\R^m}
\newcommand{\beqn}{\begin{eqnarray*}}
\newcommand{\eeqn}{\end{eqnarray*}}

\IEEEoverridecommandlockouts                              

\newcommand{\halmos}{\rule{1ex}{1.4ex}}
\newenvironment{myproof}{\noindent {\em Proof}.\ }{\hspace*{\fill}$\halmos$\medskip}
\newcommand{\epr}{\end{myproof}}
\newcommand{\bpr}{\begin{myproof}}

\title{\LARGE \bf
Global exponential stability (and contraction of an unforced system) does not imply entrainment to periodic inputs
}

\author{Alon Duvall$^{1}$ and Eduardo Sontag$^{2}$
\thanks{This work was partially supported by grants
AFOSR FA9550-21-1-0289 and NSF/DMS-2052455}%
\thanks{$^{1}$Northeastern University
        {\tt\small duvall.a@northeastern.edu}}%
\thanks{$^{2}$Northeastern University
        {\tt\small e.sontag@northeastern.edu, sontag@sontaglab.org}}%
}

\begin{document}

\maketitle
\thispagestyle{empty}
\pagestyle{empty}

\begin{abstract}

It is often of interest to know which systems will approach a periodic trajectory when given a periodic input. Results are available for certain classes of systems, such as contracting systems, showing that they always entrain to periodic inputs. In contrast to this, we demonstrate that there exist systems which are globally exponentially stable yet do not entrain to a periodic input. This could be seen as surprising, as it is known that globally exponentially stable systems are in fact contracting with respect to some Riemannian metric.

\end{abstract}

\section{INTRODUCTION}

Entrainment to a periodic input can be roughly described (a precise definition is given below) as the property that a time-independent system will converge to a periodic trajectory with the same period as a given forcing periodic input. When this occurs is often of interest in physical systems, see, for instance, \cite{article}. Of course, this property is true for stable linear systems. More generally, it is known that if a system is contractive with respect to a logarithmic norm, then it must entrain to periodic inputs, and a similar result holds for systems that are contractive with respect to arbitrary Riemannian structures \cite{7039986} \cite{FB-CTDS}. Furthermore, it is also known that if a system is globally exponentially stable to an equilibrium (GES), then the system is contractive with respect to a suitable Riemannian metric     \cite{slotine1998,Andrieu2017}. Stated in this vague fashion, it would seem that any GES system must entrain to periodic inputs. We show by means of a counterexample that this implication is false, exhibiting a GES system and a periodic input to which the system does not entrain. 

We identify the gap in the above reasoning: GES implies contractivity with respect to an arbitrary Riemannian metric, but contractivity in the absence of an input is not equivalent to contractivity when an input is present. On the other hand, for constant metrics, defined by logarithmic norms,  contractivity of the unforced system implies contractivity of the system with inputs. We go on to show that the key property needed is uniform contractivity with respect to any constant input.

\section{Background and notation}

First, we will rigorously define entrainment:

\begin{definition} 
    We say a function $v: [0,\infty)  
    \rightarrow \mathbb{R}^p$ is periodic with period $T>0$ if $v(t) = v(t+T)$ for all $t\geq0$.
\end{definition}

We will apply the above definition both to inputs ($p=m$ below) and states ($p=n$).

\begin{definition}
    Consider a system $\dot{x} = f(x,u)$, with $f: \Rn\times\Rm \to \Rn$. We say that this system \textit{entrains to periodic inputs} if the following property holds: given a function (an ``input'' or ``control'') $u: [0,\infty) \to \R^m$ which is periodic with period $T$, all solutions of $\dot{x}(t) = f(x(t),u(t))$ converge to a unique limit cycle with period $T$.
\end{definition}

We assume that the system dynamics $f$ satisfies conditions for existence and uniqueness of solutions, and forcing functions are measurable essentially bounded, see e.g.\  \cite{mct}.

We next define globally exponentially stable systems:

\begin{definition}
    Consider a system $\dot{x} = f(x,t)$, with $f: \Rn\times\R \to \Rn$ and $f(0,t)=0$ for all $t$. We say that this system is \textit{globally exponentially stable}, or just \textit{GES}, if there exists a number $\lambda > 0$ such that, for any trajectory $x(t)$ and any time $t_0$, 
    \[
    \|x(t)\| \leq e^{-\lambda (t-t_0)} \|x(t_0)\| \; \mbox{for all}\, t\geq t_0.
    \]
\end{definition}

Here $\|.\|$ is the usual Euclidean norm. We use Euclidean norm only for simplicity, since our purpose is to construct counterexamples, but more general norms could be used as well, or simply distances from $x(t)$ to $0$ in an arbitrary metric space. In our examples, the vector field $f$ will be independent of $t$, in which case we only need to consider $t_0=0$.

\section{Losing entrainment and contraction}

\subsection{Example of a GES system which does not entrain}

Consider the following two-dimensional system with two-dimensional input ($n=m=2$):
\beqn
\dot{x} &= &-x + \frac{x}{2} \sin(x^2 + y^2) -y \;+\; u_1(t)\\
\dot{y} &=& -y + \frac{y}{2} \sin(x^2 + y^2) + x  \;+\; u_2(t)\,.
\eeqn

Let $r^*$ be any positive value of $r$ at which the function
\[
f(r) = - r + \frac{r}{2}\sin(r^2)
\]
attains a strict local maximum.
In other words, $r^*$ must satisfy 
\[
\frac{\sin\left(r^2\right)}{2}+r^2\,\cos\left(r^2\right)=1
\]
and
\[
3r\cos\left(r^2\right)-2r^2\,\sin\left(r^2\right) < 0\,.
\]
The smallest local maximum is found numerically to be approximately
\[
r^* = 2.79098840365914 \,.
\]
Consider the periodic control (of period $2\pi$)
\[
u(t) = \begin{pmatrix}
    u_1(t) \\
    u_2(t)
\end{pmatrix} = -\begin{pmatrix}
    (- r^* + \frac{r^*}{2}\sin((r^*)^2)) \cos(t) \\
    (- r^* + \frac{r^*}{2}\sin((r^*)^2)) \sin(t)
\end{pmatrix}\,.
\]
The resulting system, once that this input is plugged in, is as follows:
\begin{eqnarray}
\dot{x} &=& -x + \frac{x}{2} \sin(x^2 + y^2) -y \nonumber \\
&& - (- r^* + \frac{r^*}{2}\sin((r^*)^2)) \cos(t) \label{equ:system equations 1} \\
\dot{y} &=& -y + \frac{y}{2} \sin(x^2 + y^2) + x \nonumber \\
&& - (- r^* + \frac{r^*}{2}\sin((r^*)^2)) \sin(t) \label{equ:system equations 2}
\end{eqnarray}

We show next that the curve $\gamma(t) = (r^*\cos(t),r^*\sin(t))$ is a periodic trajectory (of period $2 \pi$) of our system with the given input. Indeed, plugging in the values $(r^*\cos(t),r^*\sin(t))$  we get:
\beqn
    \dot{x} \!\! &=& \!\!\!\!
    -r^*\cos(t) + \frac{r^*\cos(t)}{2} \sin((r^*\cos(t))^2 + (r^*\sin(t))^2)\\ 
    && \quad -r^*\sin(t) - (- r^* + \frac{r^*}{2}\sin((r^*)^2)) \cos(t) \\
    &=& -r^* \sin(t)\\
\dot{y} \!\! &=& \!\!\!\!
-r^*\sin(t) + \frac{r^*\sin(t)}{2} \sin((r^*\cos(t))^2 + (r^*\sin(t))^2) \\
    && \quad + r^*\cos(t) - (- r^* + \frac{r^*}{2}\sin((r^*)^2)) \sin(t)\\
    &=& r^* \cos(t)\,.
\eeqn
Note that the right hand side of the above two equations is just the tangent vector to $\gamma(t)$, and thus we can conclude that $\gamma(t)$ is in fact a trajectory. Its image is the circle of radius $r^*$ centered at the origin.

\begin{figure}[ht]
    \centering
    \includegraphics[width=0.5\textwidth]{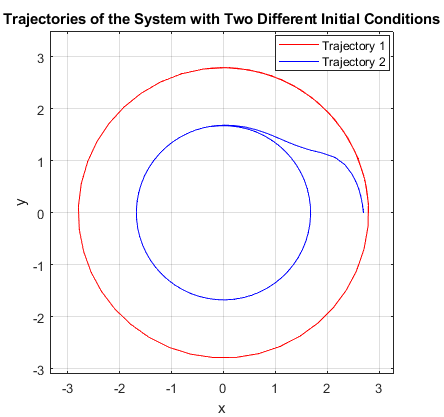}
    \caption{We plot two trajectories for our system defined by equations \ref{equ:system equations 1} and \ref{equ:system equations 2}: One with initial conditions at $(r^*,0)$ and one with initial conditions at $(r^*-0.1,0)$. We see the trajectory starting at $(r^*-0.1,0)$ does not approach, and in fact diverges from, the periodic orbit corresponding to the circle of radius $r^*$.}
    \label{fig:non_entraining}
\end{figure}

Next, let us consider an arbitrary trajectory that starts at a point of the form $(x,y)$ such that $x^2 + y^2 <(r^*)^2$. We will use polar coordinates ($r(t),\theta(t))$ to represent this trajectory, so that
$x(t) = r(t) \cos(\theta(t))$
and $y(t) = r(t) \sin(\theta(t))$.
Thus $r(t)=\sqrt{x(t)^2+y(t)^2}$ along this trajectory and the initial condition $(x(0),y(0))$ has the form in polar coordinates
$(r(0),\theta(0))$, with $ r^2(0) = x(0)^2 + y(0)^2$.
We will show that $r(t)$ is decreasing whenever $r(t) < r^*$ is sufficiently close to $r^*$. We have that (note we are now suppressing that $x,y$ and $r$ are all functions of $t$):

\begin{align*}
\frac{\dot{r^2}}{2} &= x \dot{x} + y\dot{y} \\
&= -(x^2 + y^2)  + \frac{x^2 + y^2}{2} \sin(x^2 + y^2) \\
&\phantom{=} - (x \cos(t) + y \sin(t))(- r^* + \frac{r^*}{2}\sin((r^*)^2)) \\
&= -r^2 + \frac{r^2}{2} \sin(r^2) \\
&\phantom{=} - (r \cos(t - \theta)) (- r^* + \frac{r^*}{2}\sin((r^*)^2)).
\end{align*}
If $r$ is close enough to $r^*$ such that $0 > f(r) - f(r^*) $, then
\begin{align*}
&-(r^2) + \frac{r^2}{2} \sin(r^2) - (r \cos(t - \theta)) (- r^* + \frac{r^*}{2}\sin((r^*)^2)) \\
&\leq -(r^2) + \frac{r^2}{2} \sin(r^2) - (r) (- r^* + \frac{r^*}{2}\sin((r^*)^2)).
\end{align*}
\[
= r(f(r) - f(r^*))< 0 \,.
\]
Thus we see that any points in the interior of the circle $r = r^*$ and close enough to the trajectory $\gamma(t)$ actually move away from from $\gamma(t)$, and thus we do not approach $\gamma(t)$. Figure~\ref{fig:non_entraining} illustrates this.

This system is equivalent to the following system in polar coordinates:
\beqn
\dot{r} &=& f(r) \;=\; - r + \frac{r}{2}\sin(r^2)\\
\dot{\theta} &=& 1\,.
\eeqn
Indeed, the motivation for this system was that we simply want to ``force'' an unstable periodic orbit, such that points slightly closer to the origin move away from this periodic orbit, showing that there is no global convergence to a unique periodic orbit. This form of the system also makes it easy to see that it is in fact GES, since we have that 
$f(r) \leq -\frac{r}{2}$ for all nonnegative $r$, which implies that, for all solutions,
\[
\|x(t)\| \leq e^{-\frac{1}{2} t } \|x(0)\| \;\; \mbox{for all} \; t\geq 0\,.
\]

\subsection{Example of losing contraction upon translation}

A definition of contraction with respect to a state-dependent (and time-independent) Riemannian metric on $\mathbb{R}^n$ is given in \cite{slotine1998}. We prefer to use the form given in \cite{parillo_slotine_2008}, Definition 1:

Given the $n$-dimensional autonomous system 
\[
\dot{x} \;=\; f(x) \,,
\] 
a contraction metric is an $n \times n$ symmetric matrix $M(x)$ that is uniformly positive definite (that is, 
$v^T M(x) v \geq a$  for all $x$ and $v$ in $R^n$, 
for some $a>0$) such that
\[
\frac{\partial f}{\partial x}^T M(x) + M(x) \frac{\partial f}{\partial x} + \dot M (x)
\]
is uniformly negative definite, where the notation $\dot M (x)$ is shorthand for the matrix whose $(i,j)$th entry is
\[
\dot{M}(x)_{ij} = \frac{\partial M_{ij}}{\partial x}^T f(x) \,.
\]
Theorem 1 in \cite{parillo_slotine_2008} says that if a contraction matrix exists and there is some equilibrium, then all trajectories converge to it.
Moreover, under the stronger assumption that 
\[
\frac{\partial f}{\partial x}^T M(x) + M(x) \frac{\partial f}{\partial x} + \dot M (x) \leq \;-\; \alpha M(x)
\]
for all $x$, for some $\alpha>0$ (where $A\leq B$ means that $B-A$ is positive definite) then there is automatically an equilibrium.
The second result follows from the fact that $d(x(t),y(t)) \leq e^{- (\alpha/2)t} d(x(0),y(0))$ along trajectories, where $d$ is the geodesic distance associated to the metric $M(x)$.

Let us consider now as an example the scalar system $\dot{x} = f(x) = - x + \frac{x}{2}\sin(x^2)$ where  $n=1$. One needs a function $m:\R \rightarrow \R$ such that $m(x)\geq a > 0$ for all $x$ and 
\[
m'(x)f(x) + 2f'(x)m(x) \;\leq\; - \beta m(x)
\]
for all $x$, where $\beta> 0$.
In our case, we pick 
$m(x) = 1/(\sin(x^2)/2 - 1)^2$,
which is larger than 4/9 for all $x$, and with
$f(x) = (x \sin(x^2))/2 - x$
we have that 
\beqn
m'(x)f(x) + 2f'(x)m(x) &=& \frac{4}{\sin(x^2) - 2}\\
&\leq&  \frac{-1/3}{(\sin(x^2)/2 - 1)^2} \,.
\eeqn
Thus we see that for $\beta = 1/3$ our inequality is satisfied. Indeed, this follows since 
\[
\frac{4(\sin(x^2)/2 - 1)^2}{\sin(x^2) - 2} \leq \frac{1}{\sin(x^2) - 2} \leq -\frac{1}{3} \,.
\]

We next show that contractivity breaks down for this metric if we add a constant input, i.e., we will consider the system $\dot{x} = f(x) + c$. In this case we have that
\beqn
m'(x) (f(x)+c) + 2f'(x) m(x) &=& \frac{4}{\sin(x^2) - 2} \\
&&\quad + c\frac{16 x \cos(x^2)}{(2 -\sin(x^2))^3} \\
&\geq& c\frac{16 x \cos(x^2)}{(2 -\sin(x^2))^3} - 4 \\
&\geq& c\frac{16 x \cos(x^2)}{27} - 4\,.
\eeqn
If we pick $c = \frac{27}{16}$ and $x = 4\sqrt{2 \pi }$ we see that our final expression is clearly positive, and thus our system does not contract everywhere with respect to this metric (despite the fact that when $c=0$, our system \textit{is} contractive with respect to this metric).

\section{Constant and non-constant inputs, and entrainment}

In this section, we will first remark that uniform contractivity (and hence entrainment) under \emph{constant} (additive) inputs that take values in a given set $B$ implies uniform contractivity (and hence entrainment) under \emph{arbitrary} inputs taking values in $B$. We will then prove, however, that if $B=\R^n$ then the only metrics for which this property can hold are constant, but if $B$ is bounded, there do exist non-constant metrics that provide contractivity.

\subsection{Connecting constant and non-constant inputs}

\begin{thm}
\label{thm:additive}
    Consider the system with additive inputs $\dot{x} = f(x) + u$, where the inputs take values $u(t)\in U \subset \R^n$.
    Suppose that this system is uniformly contractive with respect to a metric $M(x)$, uniformly over all constant inputs $u(t)\equiv c\in U$. Then the system is uniformly contractive with respect to a metric $M(x)$, uniformly on all inputs $u$ with values in $U$.
\end{thm}

\bpr
The assumption means that there is a constant $\beta>0$ such that, for all $x$, and all $g(x) = f(x)+c$, $c\in U$, the following inequality holds:
\[
\frac{\partial g}{\partial x}^T M + M \frac{\partial g}{\partial x} + \dot{M} \leq -\beta M
\]
(where $M = M(x)$).
Here
\[
\dot{M} = \dot{M}_1 + \dot{M}^{[c]}_2\,,
\]
where we define 
\[
\{\dot{M}_1\}_{ij} := \frac{\partial M_{ij}}{\partial x} f(x),\, \{\dot{M}^{[c]}_2\}_{ij} := \frac{\partial M_{ij}}{\partial x} c \,.
\]
Thus we have an estimate, uniform on $x$ and $c$,
\[ 
\frac{\partial f}{\partial x}^T M + M \frac{\partial f}{\partial x} + \dot{M}_1 + \dot{M}^{[c]}_2 \leq -\beta M.
\]
Now given an arbitrary input $u$ with values in $U$, we'd like to have 
\[ 
\frac{\partial f}{\partial x}^T M + M \frac{\partial f}{\partial x} + \dot{M}_1 + \dot{M}^{[u(t)]}_2 \leq -\beta M
\]
for all $t$,
where $\{\dot{M}^{[u(t)]}_2\}_{ij} := \frac{\partial M_{ij}}{\partial x} u(t)$.
This is true because $u(t)\in U$.
\epr

As a corollary we have, under the same assumptions, entrainment to any periodic input $u(t)$ satisfying $u(t) \in U$ for all $t \geq 0$.
This follows from a result from \cite{slotine1998} which states that if a system $\dot{x} = f(x,u(t))$ is contracting then we have entrainment.

In Section~\ref{sec:metric} we will provide a version of Theorem~\ref{thm:additive} which applies to non-additive inputs, and is formulated in an abstract metric space setup.

\subsection{Metrics must be constant for contractivity under arbitrary additive constant inputs}

We next show that when $U=\R^n$, uniform contraction metrics under aditive inputs must be constant.

\begin{thm}
\label{thm: all c implies constant M}
    Suppose we have a system $\dot{x} = f(x)$ which is contractive with respect to some metric $M(x)$, where $M(x)$ is not constant. Then there exists a constant input $c \in \mathbb{R}^n$ such that $\dot{x} = f(x) +c$ is no longer contractive with respect to $M(x)$.
\end{thm}

\bpr
Consider the system $\dot{x} = f(x) + c$, where $c$ is a constant vector. Since we are free to choose $c$ as we please, we will occasionally change $c$ in the course of the proof. For this system to be contractive we must have that
\[
\frac{\partial f}{\partial x}^T M + M \frac{\partial f}{\partial x} + \dot{M} \leq 0.
\]
let us write $\dot{M} = \dot{M}_1 + \dot{M}_2$ where $\{\dot{M}_1\}_{ij} = \frac{\partial M_{ij}}{\partial x} f(x)$ and $\{\dot{M}_2\}_{ij} = \frac{\partial M_{ij}}{\partial x} c$. Now we have that
\[
\frac{\partial f}{\partial x}^T M + M \frac{\partial f}{\partial x} + \dot{M} = \frac{\partial f}{\partial x}^T M + M \frac{\partial f}{\partial x} + \dot{M}_1 + \dot{M}_2.
\]
Since $M(x)$ is not constant we can always find a vector $x$ and a vector $c$ such that $\dot{M_2}(x)$ is nonzero, as well as a vector $z$ such that $ \alpha = z^t\dot{M}_2(x) z \neq 0 $. If $\alpha < 0$, then replace $c$ with $-c$ so that we will instead get that $\alpha > 0$. Now assume $\alpha > 0$ and set
\[\beta = z^t (\frac{\partial f}{\partial x}^T M + M \frac{\partial f}{\partial x} + \dot{M}_1) z \]
Pick $N$ such that $N\alpha > |\beta|$. Upon replacing $c$ with $Nc$, we have that

\[ z^t (\frac{\partial f}{\partial x}^T M + M \frac{\partial f}{\partial x} + \dot{M}_1 + \dot{M}_2) z  = \beta + N \alpha > 0.\]

Thus we see that given that $M(x)$ is nonconstant, we can always find a constant input $c$ for which our system is no longer contractive everywhere with respect to $M(x)$.
\epr

\subsection{Bounded additive controls have non-constant contractions}

One might wonder if we can do better than the previous theorem. Along these lines, one might ask: Does there exist a \emph{bounded} set $B$ such that if we are uniformly contractive for inputs taking values in $B$ with respect to a certain metric $M$, then our metric \emph{must} be constant? This is not true, as we will now show by means of a counterexample.

\begin{proposition}
    Given an arbitrary bounded set $B \subseteq \mathbb{R}$ there exists a one dimensional system $\dot{x} = f(x) + u(t)$ such that if $u(t) \in B$ for all $t \geq 0$, then there exists a non-constant metric for which this system is contractive. 
\end{proposition}

\bpr
    Using the notation from Theorem \ref{thm: all c implies constant M} we want to establish the inequality
\[
\frac{\partial f}{\partial x}^T M + M \frac{\partial f}{\partial x} + \dot{M}_1 + \dot{M}_2 \leq -\beta M
\]
for some particular choices of (non-constant) $M$ and $f$.
    
Consider the system where $f(x) = -x + c$. For this our equation looks like
\[
(c - x) \frac{\partial M}{\partial x} \leq( 2-\beta) M \,.
\]
Suppose $\beta = 1$ and $M(x) = 1 + \epsilon(x)$. Then we have that our inequality is equivalent to
\[
(c - x) \epsilon'(x) \leq 1 + \epsilon(x)\,.
\]
    Consider the quantity $\epsilon_m(x) = e^{-x^2/m}$. Substituting $\epsilon_m(x)$ for $\epsilon(x)$ in our inequality we have that
\[
(c - x) \epsilon_m'(x) = \frac{-2x(c-x)}{m} e^{-x^2/m}.
\]
Note that for large enough $m$, this expression goes uniformly to 0 (because $B$ is bounded), and so for large enough $m$ we will have that $(c - x) \epsilon_m'(x) \leq 1 + \epsilon_m(x)$. Thus we have that $M(x) = 1 + \epsilon_m(x)$ is a valid contraction metric for large enough $m$.
\epr

\section{Generalizations and comments}

We discuss now several directions in which our results can be generalized.

\subsection{A more general sufficient condition forcing a metric to be constant}

We will use the following notation:
\begin{enumerate}
    \item $\lambda_{\max}(A)$ is the maximum eigenvalue of an arbitrary matrix $A$.
    \item $\|v\|_2$ is simply the Euclidean norm of an arbitrary vector $v$.
    \item  $\|A\|_2 = \sqrt{\lambda_{\max} (A^T A)}$ for an arbitrary matrix $A$.
    \item $\mu_2(A) = \lambda_{\max}\left(\frac{A + A^t}{2} \right)$ for an arbitrary matrix $A$.
\end{enumerate}

The last three operations are a vector norm, its induced matrix norm, and its induced logarithmic norm, respectively.
Logarithmic norms are routinely used in contraction theory, and an early exposition is in the control textbook~\cite{DesoerVidyasagar}.
The following elementary facts are well-known, but we include proofs to make the exposition self-contained.

\begin{lem}
    If $A$ and $B$ are symmetric and $A \leq B$, then $\mu_2(A) \leq \mu_2(B)$.
\end{lem}

\bpr
    First note $\mu_2(A)$ is simply the maximum value of its eigenvalues, for a symmetric matrix $A$. Suppose $A$'s maximum eigenvalue is larger than $B$'s maximum eigenvalue. Then we have that if $x$ is the corresponding eigenvector for $A$ that
\beqn
\frac{x^t (B- A) x}{x^t x} &=& \frac{x^t B x}{x^t x} - \lambda_{\max}(A)\\
&\leq& \lambda_{\max}(B) - \lambda_{\max}(A)< 0.
\eeqn
This a contradiction, and so we are done.
\epr

\begin{lem}
    We have that $\mu_2(AB) \leq \|A \|_2 \|B \|_2$.
\end{lem}

\bpr
    Note that for a symmetric matrix $S$ we have that $\lambda_{\max}(S) \leq |\lambda_{\max}(S)| \leq \|S\|_2$. We have that
    \beqn
        \mu_2(AB) &=& \lambda_{\max} \left(\frac{AB + (AB)^t}{2}\right) \leq \left\| \frac{AB + (AB)^t}{2} \right\|_2 \\
        &\leq& \|A \|_2 \|B \|_2
    \eeqn
\epr

    Recall that the \textit{convex closure} of a set of point $P$ is indicated by $\mbox{conv}(P)$ and is the set of all points $p$ such that we can find a finite set of points $O \subset P$ such that $p = \sum_{i} \lambda_i o_i $, where $\lambda_i \geq 0, \sum_i \lambda_i = 1$ and $o_i \in O$.

\begin{thm}
\label{thm:general constant metric}

    Suppose we have: 

    \begin{enumerate}
        \item An $n$ dimensional system $\dot{x} = f(x,u)$.
        \item $k$ infinite sequences of constant inputs $\{u_{i,j}\}$ where $1 \leq i$ and $1 \leq j \leq k$.
        \item For each fixed $x$ there exists a large enough $i_0$ so that $\mbox{conv}(\{\frac{f(x,u_{i,j})}{\|f(x,u_{i,j})\|_2}\}_{1 \leq j \leq k})$ always contains a certain open sphere $O$ centered at 0 for $i \geq i_0$.
        \item For any fixed $j$ and for all $x$ we have that.
    \[\lim_{i \rightarrow \infty} \frac{\left\| \frac{\partial f(x,u_{i,j})}{\partial x}\right\|_2 }{\|f(x,u_{i,j})\|_2}  = 0.\]
    \end{enumerate}

    Then if our system is contractive with respect to a Riemannian metric for each of our constant inputs, the metric must be constant.
\end{thm}

\bpr
    First we fix an arbitrary $x \in \mathbb{R}^n$. We will carry out all our computations at this set $x$. We have our necessary equation for contractivity:
    \[
\frac{\partial f}{\partial x}^T M + M \frac{\partial f}{\partial x} + \dot{M} \leq 0.
\]
Now we can rearrange to get
\[
 \dot{M} \leq -\frac{\partial f}{\partial x}^TM  - M \frac{\partial f}{\partial x}  .
\]
Now upon taking the logarithmic norm $\mu_2$ of both sides we have that
\beqn
    \mu( \dot{M} )_2 &\leq & \mu_2 \left(-\frac{\partial f}{\partial x}^TM  - M \frac{\partial f}{\partial x}  \right) \\
    &\leq &  \mu_2 \left(-\frac{\partial f}{\partial x}^TM \right) + \mu_2\left(- M \frac{\partial f}{\partial x}\right) \\
    &\leq & 2\left\|-\frac{\partial f}{\partial x}^T \right\|_2 \|M\|_2 \\
    &=& 2\left\|\frac{\partial f}{\partial x}^T \right\|_2 \|M\|_2.
\eeqn
Dividing by $\|f(x,u_{i,j})\|_2$ we get
\[
\mu_2\left(\frac{ \dot{M} }{\|f(x,u_{i,j})) \|_2}\right)  \leq \frac{2 \|M\|_2 \| \frac{\partial f}{\partial x}^T \|_2  }{\|f(x,u_{i,j})\|_2} \xrightarrow[i \rightarrow \infty]{} 0.
\]
Thus we see that we have $\lim_{i \rightarrow \infty} \mu_2\left(\frac{ \dot{M} }{\|f(x,u_{i,j})\|_2} \right) \leq 0$. Now the $lk$'th entry of $\frac{ \dot{M} }{\|f(x,u_{i,j})\|_2}$ is $\frac{\partial M_{lk}}{\partial x} \frac{f(x,u_{i,j})}{\|f(x,u_{i,j})\|_2}$. Here $\frac{\partial M_{lk}}{\partial x}$ is the gradient of $M_{lk}$ and is a row vector with $n$ entries. Suppose we have a matrix $G$ such that its $lk$'th entry is $\{G\}_{lk} = \frac{\partial M_{lk}}{\partial x} v $ where $v \in \mathbb{R}^n$. Let $\alpha$ be positive real number such that $v \in \alpha O$. Now we can write
\[ v = \sum_j \alpha_{i,j} \frac{f(x,u_{i,j})}{\|f(x,u_{i,j} \|}, \]
where we have that $0 \leq \alpha_{i,j}  \leq \alpha$ for large enough $i$ and all $j$ (this follows from the fact that $v \in \alpha O \in \alpha* \mbox{conv}(\{\frac{f(x,u_{i,j})}{\|f(x,u_{i,j})\|_2}\}_{1 \leq j \leq k})$). Thus we have that
\[\{G\}_{lk} = \frac{\partial M_{lk}}{\partial x} v  =  \sum_j \alpha_{i,j} \frac{\partial M_{lk}}{\partial x} \frac{f(x,u_{i,j})}{\|f(x,u_{i,j})\|_2}. \]
Now let $M^{i,j}$ be the matrix with $\{M^{i,j}\}_{lk} =  \frac{\partial M_{lk}}{\partial x} \frac{f(x,u_{i,j})}{\|f(x,u_{i,j})\|_2}$. We have that
\[\mu_2(G) \leq \alpha \sum_j \mu_2(M^{i,j}) \xrightarrow[i \rightarrow \infty]{} 0.\]
Thus we have that $\mu_2(G) \leq 0$ and by replacing $v$ with $-v$ we also have $\mu_2(-G) \leq 0$ which only happens when $G = 0$. Thus we have that $\frac{\partial M_{lk}}{\partial x} v = 0$ for all $v$, and so we must have $\frac{\partial M_{lk}}{\partial x} = 0$. Since this is true at all $x$ we must have that $M(x)$ is the constant matrix.

\subsubsection{Example}

Consider the system
\[\begin{bmatrix}
    \dot{x} \\
    \dot{y} \\
    \dot{z} 
\end{bmatrix} = \begin{bmatrix}
    -1 & 0 & 0 \\
    0 & -1 & 0 \\
    0 & 0 & -1 \\
\end{bmatrix} \begin{bmatrix}
    x \\
    y \\
    z
\end{bmatrix} + \begin{bmatrix}
    1 & 0 & 0 & -1 \\
    0 & 1 & 0 & -1 \\
    0 & 0 & 1 & -1
\end{bmatrix} \begin{bmatrix}
    u_1 \\
    u_2 \\
    u_3 \\
    u_4
\end{bmatrix}.\]
Suppose we have that for all $k$ that $u_k \geq 0$, and at most one of $u_1,u_2,u_3$ or $u_4$ are nonzero. Take as four sequences of controls 
\[\{u_{k,1}, u_{k,2}, u_{k,3}, u_{k,4} \} = \left\{\begin{bmatrix} k \\ 0 \\ 0 \\ 0 \end{bmatrix},\begin{bmatrix} 0 \\ k \\ 0 \\ 0 \end{bmatrix},\begin{bmatrix} 0 \\ 0 \\ k \\ 0 \end{bmatrix},\begin{bmatrix} 0 \\ 0 \\ 0 \\ k \end{bmatrix} \right\}.\]
Where $k$ is any positive integer. This sequence of controls satisfies the conditions of our theorem. Indeed, we have that $\frac{f(x,u_{k,j})}{\|f(x,u_{k,j})\|_2}$ is simply a constant vector for each $j$, and thus the convex closure of all 4 of these constant vectors always contains an open ball. We also have that
\[\lim_{i \rightarrow \infty} \frac{\left\| \frac{\partial f(x,u_{k,j})}{\partial x}\right\|_2 }{\|f(x,u_{k,j})\|_2}  = 0\]
since the numerator is always a constant, while the denominator always goes to infinity.

Thus while we are contractive with regard to a constant metric (such as the usual Euclidean metric), we are not contractive with respect to a nonconstant metric.

\subsubsection{Example}

Consider the system 
\[\dot{x} = f(x) + c.\]
Here $c \in \mathbb{R}^n$ is our constant control, which we allow to take on any value. We see that we have
    \[\lim_{\|c\|_2 \rightarrow \infty} \frac{\left\| \frac{\partial f(x,c)}{\partial x}\right\|_2 }{\|f(x,c)\|_2} = \lim_{\|c\|_2 \rightarrow \infty} \frac{\left\| \frac{\partial f(x)}{\partial x}\right\|_2 }{\|f(x)+c\|_2}  = 0\]
since the numerator does not depend on $c$, and the denominator goes to infinity as $\|c\|_2$ tends to infinity. We can also always find $n+1$ sequences of constant inputs $c_{i,j}$ such that $\lim_{i \rightarrow \infty} \|c_{i,j}\|_2 = \infty$ and $\mbox{conv}(\{\frac{f(x) + c_{i,j}}{\|f(x) + c_{i,j}\|_2} \}_{1 \leq j \leq n+1 })$ always contains a fixed open ball (e.g., take $c_{i,j} = i v_j$, where $\mbox{conv}(\{v_j\}_{1 \leq j \leq n+1})$ contains an open ball). Thus the conditions of our theorem are satisfied, and thus we have another proof of Theorem $\ref{thm: all c implies constant M}$.
\epr

\subsection{Contraction on a metric space}
\label{sec:metric}

Suppose we are given a metric space $\mathbb{M}$ with metric $d$. Refer to the space of mappings of $\mathbb{M}$ to itself as $C(M)$. 

\begin{definition}
    The \textit{topology induced by its supremum distance} on $C(M)$ is the topology induced by the metric $\mathcal{D}$ on $C(M)$ defined by
    \[\mathcal{D}(f,g) = \sup\{d(f(x),g(x)| x \in M\}.\]
\end{definition}

\begin{definition}
    Given two functions $u,v: \mathbb{R} \rightarrow \mathbb{R}^n$ their \textit{concatenation} $uv$ at $t_0$ is the function defined as\[uv = 
    \begin{cases}
        u & \text{when } t < t_0 \\
        v & \text{when } t \geq t_0
    \end{cases} \]
\end{definition}

\begin{thm}
    Suppose we have:

    \begin{enumerate}
        \item A metric space $\mathbb{M}$ with metric $d$.
        \item The space of mappings of $\mathbb{M}$ to itself $C(\mathbb{M})$ (give this space the topology induced by its supremum distance).
        \item The space of measurable, bounded, and locally integrable functions from $\mathbb{R} \rightarrow \mathbb{R}^n$ (call this space $L_M$, give this space its supremum norm, where the norm is Euclidean).
        \item A mapping $ \phi:L_M \times \mathbb{R} \times \mathbb{R} \rightarrow C(\mathbb{M})$ satisfying that 
        \begin{enumerate}
            \item It is continuous in $L_M$ (given two arbitrary and fixed entries for the other arguments)
            \item  $\phi(g,t_2,t_3) \circ \phi(f,t_1,t_2) = \phi(gf,t_1,t_3)$ where $gf$ is the concatenation of the two functions at $t_2$.
            \item $\phi(f(t),t_1,t_2) = \phi(f(t+T), t_1 -T, t_2 - T )$ for all $T,t,t_1,t_2 \in \mathbb{R}$.
            \item Whenever we have a compact set $U \subseteq L_M$ of constant mappings then there exists $\lambda < 0$ such that $d(\phi(f,0,t)(x),\phi(f,0,t)(y)) \leq e^{\lambda t}d(x,y)$ for all $x,y \in \mathbb{M}$, all $t\geq 0$, and all $f \in U$.
        \end{enumerate}

    \end{enumerate}

    Let $L_{ct}$ be the space of mappings produced by concatenating finitely many times functions from $U$, and let $\overline{L}_{ct}$ be the closure of this set.  If we have $f \in \overline{L}_{ct}$ and $t_1 > t_2$ then $\phi(f,t_1,t_2)$ is a contraction.
    
\end{thm}

\bpr
    First we will argue that all piecewise constant functions also give us contractions. This is clear from the requirement that $\phi(g,t_2,t_3) \circ \phi(f,t_1,t_2) = \phi(gf,t_1,t_3)$. Indeed, letting $\phi_g = \phi(g,t_2,t_3)$ and $ \phi_f = \phi(f,t_1,t_2)$ we have if $d(\phi_f(x),\phi_f(y)) \leq \lambda_1 d(x,y)$ and $d(\phi_g(x),\phi_g(y)) \leq \lambda_2 d(x,y)$ then we have
    \[ d(\phi_g (\phi_f(x)),\phi_g(\phi_f(y))) \leq \lambda_2 \lambda_1 d(x,y).
    \]
    Thus the composition of $\phi_f$ and $\phi_g$ still gives us a contraction and so $\phi(gf,t_1,t_3)$ must be a contraction. It follows that if we concatenation finitely many functions into a piecewise constant function $F$, we will still have that $\phi(F,t_1,t_2)$ is a contraction. Since any piecewise constant function (taking on finitely many different values) can be produced in this manner, we can conclude that for all $F \in L_{ct}$ that $\phi(F,t_1,t_2)$ is in fact a contraction.

     Now suppose we consider a compact set $B \subseteq U$ of constant functions, taking their image from a bounded codomain in $\mathbb{R}^n$. For each $f \in B$ we have that $d(\phi(f,t_1,t_2)(x),\phi(f,t_1,t_2)(x)) \leq e^{\lambda (t_2-t_1)} d(x,y)$ where $\lambda$ is the contraction constant we know will work for all our constant functions by condition 4(d).

     Now suppose we have a piecewise constant function $f$ on the interval $[t_1,t_2]$ and set $\phi(f(x),t_1, t_2) = \phi_f$. Suppose our function takes on value $c_i$ on interval $i$ of length $\alpha_i(t_2-t_1)$, where $\sum_i \alpha_i = 1$ and the $c_i$ are all contained in some compact set $B$. Composing all these pieces as in property 4(b) of our conditions, we have that
    \beqn
        d(\phi_f(x),\phi_f(y)) &\leq& (\prod_i e^{\lambda \alpha_i(t_2 -t_1)}) d(x,y)\\
        &=& e^{\lambda (t_2 - t_1)} d(x,y).
    \eeqn
    Thus not only are our piecewise constant mappings contractions, but there is a contraction constant $\lambda$ that they all satisfy (assuming we are taking our constant functions from a bounded set). 
    
    Now for a general function $f$ that can be approximated (in supremum norm) by piecewise constant functions, note that due to $f$ being bounded its values are contained in a compact set. Thus it can be approximated by piecewise constant functions taking on values from a compact set. Thus all these approximating functions have a universal contraction coefficient $\lambda < 0$, and by continuity of $\phi$ we have that $\phi(f,t_1,t_2)$ must also be a contraction. 

    Indeed, if we have a sequence of piecewise constant $f_i$ converging to $f$ so that $\phi_{f_i}$ converges to $\phi_f$, then we have that for a given pair of points $x,y \in \mathbb{M}$ and for some $\epsilon > 0$ that
    \beqn
        d(\phi_f(x),\phi_f(y))
        &\leq& d(\phi_{f}(y),\phi_{f_i}(y)) + d(\phi_{f_i}(x),\phi_{f_i}(y)) \\ && + \; d(\phi_{f_i}(x),\phi_{f}(x)) \\
        &\leq& e^{-\lambda(t_1 - t_2)} d(x,y) + \epsilon.
    \eeqn
    Thus since we can pick $\epsilon$ to be arbitrarily small (it becomes small as $i \rightarrow \infty$) we have that 
    \[d(\phi_f(x),\phi_f(y)) \leq e^{-\lambda(t_1 - t_2)} d(x,y).
    \]
\epr

This theorem, while abstract, applies concretely to ordinary differential equations. If we have a system $\dot{x} = f(x) + u(t)$ where $u(t)$ is in a compact set $B$, $f(x)$ is smooth, and for $u(t)$ constant our system is contractive, then by Theorem 55 in \cite{mct} and Theorem 2 in \cite{slotine1998} we have that $\phi$ satisfies the required conditions and so we can use our theorem. Thus if $\dot{x} = f(x) + c$ is a contraction for $c$ in some compact set $B$ and if $u(t) \in B$ for all $t \geq 0$ then our system will still be a contraction. Thus if $u(t)$ is periodic we will have entrainment.

\section{Conclusion and future directions}
 In this conference paper we studied the connections between global exponential stability, contractions with respect to constant and nonconstant metrics, and entrainment. In the full version of this paper, we will describe additional sufficient conditions for which a globally exponentially stable system would in fact entrain to periodic inputs.

\bibliographystyle{unsrt}
\bibliography{mybib}

\end{document}